\documentclass{article}
\usepackage{spconf,amsmath,graphicx,cite}
\usepackage{amsfonts,amssymb}

\title{Relevance Singular Vector Machine for Low-rank Matrix Sensing}
\name{Martin Sundin, Saikat Chatterjee, Magnus Jansson and Cristian R. Rojas\thanks{This work was partially supported by the Swedish Research Council under contract 621-2011-5847.}}
\address{ACCESS Linnaeus Center, School of Electrical Engineering\\
KTH Royal Institute of Technology, Stockholm, Sweden\\
{\small \tt masundi@kth.se, sach@kth.se, janssonm@kth.se, crro@kth.se}}

\begin{document}

\maketitle
\begin{abstract}
In this paper we develop a new Bayesian inference method for low rank matrix reconstruction. We call the new method the Relevance Singular Vector Machine (RSVM) where appropriate priors are defined on the singular vectors of the underlying matrix to promote low rank. To accelerate computations, a numerically
efficient approximation is developed. The proposed algorithms are applied to matrix completion and matrix reconstruction problems and their performance is studied numerically.
\end{abstract}
\section{Introduction}
\label{sec:intro}

Low-rank matrix reconstruction (LRMR) from under-sampled measurements has attracted considerable interest
in recent times \cite{Candes09,Candes10,Fazel,Babacan}. An example of LRMR is matrix completion \cite{Candes09,Candes10}. 
LRMR can be viewed as a further generalization of the thriving topic called compressed sensing \cite{Baraniuk}.
Typically reconstruction algorithms are of two main types, convex and greedy methods. 
Convex methods try to minimize the nuclear norm of the underlying matrix \cite{Candes09,Candes10}.
Under some conditions, the convex relaxation algorithms are optimal. 
On the other hand 
greedy search techniques also have been developed based on heuristics, providing reasonable practical performance.

In estimation theory, the Bayesian framework has a strong role. While attempts have been made 
to develop non-Bayesian tools for LRMR -- such as convex and iterative greedy techniques -- not 
much effort has been made to develop a fully Bayesian framework.
A recent attempt 
in this direction
is \cite{Babacan}. In \cite{Babacan} an indirect 
approach is 
based on 
reparameterization of the original low-rank matrix 
by two component matrix factorization and then 
inducing low-rank structures on the two matrices 
through standard Variational Bayesian tools. 
Hence, the approach of \cite{Babacan} does not use a low-rank prior for the matrix itself.
The method of \cite{Babacan} provides good performance when the matrix is nearly square, but suffers in 
performance for skewed dimensional matrices.

In this paper we attempt to use a low-rank prior directly for the matrix itself. 
The low-rank prior is constituted via the property that the number of relevant singular vectors (singular vectors spanning the range space of the matrix) 
are limited for low-rank matrices. 
We describe the singular vectors by introducing precisions (inverse variances) for the left and right singular vectors. 
Since our method naturally generalizes the Relevance Vector Machine (RVM) for sparse vectors \cite{Tipping,Bishop} to low-rank matrices we refer to it as the Relevance Singular Vector Machine (RSVM). 
The development of the RSVM comes with a non-trivial 
prior formulation followed by the MAP rule for parameter estimation. 
Through several experiments we compare the efficiency of the RSVM to the Variational Bayesian based
technique of \cite{Babacan} and the standard nuclear norm minimization technique.

%
%

\section{Problem formulation}
\label{sec:format}

In our model a low rank matrix $\mathbf{X} \in\mathbb{R}^{p\times q}$ is observed through noisy linear measurements
\begin{align}
\label{measurements}
\mathbf{y} =  \mathbf{A}\mathrm{vec}(\mathbf{X}) +  \mathbf{n}
\end{align}
where $\mathbf{y} \in \mathbb{R}^m$ and $\mathbf{A} \in \mathbb{R}^{m \times pq}$. The rank of the matrix is low, $\mathrm{rank}(\mathbf{X}) = r \ll \min(p,q)$, and unknown. The noise is zero-mean white Gaussian, $E[\mathbf{nn^\top}] = \sigma_n^2 \mathbf{I}_m$, with unknown variance $\sigma_n^2$. In the matrix completion scenario \cite{Candes10}, the linear operator $\mathbf{A}$ chooses a set of elements from $\mathbf{X}$. The goal is to reconstruct the matrix $\mathbf{X}$ from the measurements $\mathbf{y}$.

\subsection{Prior work}

The problem of reconstructing $\mathbf{X}$ from \eqref{measurements} has been studied using several approaches \cite{Candes09,Candes10,Fazel,Babacan,Zachariah}. One convex approach is the relaxed nuclear norm heuristic \cite{Candes10}, which sets
\begin{align}
\label{nuclearnorm}
\mathbf{\hat{X}} = & \arg \min_{\mathbf{X}} ||\mathbf{X}||_* , \,\,  \text{ such that } ||\mathbf{y} - \mathbf{A}\mathrm{vec}(\mathbf{X})||_2 \leq \delta 
\end{align}
where $\delta \geq 0$ is a regularization parameter and $||\mathbf{X}||_*$ is the nuclear norm \cite{Candes09}. The estimate \eqref{nuclearnorm} depends on the parameter $\delta$, which needs to be set in advance. Greedy methods for \eqref{measurements} have also been developed \cite{admira,Zachariah}. However, most greedy methods require either the rank or noise variance to be known a priori.


LRMR can be viewed as a further generalization of compressed sensing. In Bayesian Compressed Sensing (BCS) \cite{bayesianCS}, the Relevance Vector Machine (RVM) \cite{Tipping} is used to estimate a sparse vector from an underdetermined set of measurements. The RVM uses hyperpriors (or hierarchical priors) to promote sparsity. That is, the distributions of the hyperpriors promote sparsity in the estimate. The RVM does not require the sparsity or noise variance to be known a priori.

One analogue of BCS for matrix completion and robust PCA was developed by Babacan et. al. \cite{Babacan}. In \cite{Babacan}, the matrix $\mathbf{X}$ was factorized as
\begin{align}
\label{factorization}
\mathbf{X} = \mathbf{LR^\top}
\end{align}
where $\mathbf{L} \in \mathbb{R}^{p \times s}$, $\mathbf{R} \in \mathbb{R}^{q \times s}$ and $s$ is a user defined parameter. The authors used the variational Bayesian framework to iteratively estimate the column vectors of $\mathbf{L}$ and $\mathbf{R}$. In this model, low-rank is promoted by reparameterization the problem and using sparsity inducing priors, rather than using low-rank inducing priors for the matrix itself.




\section{Relevance Singular Vector Machine} 
\label{sec:typestyle}



To generalize BCS to low-rank matrix reconstruction, we construct matrix precisions that induce low-rank. The precisions give information about the low-rank structure of the matrix, i.e., the left and right singular vectors of the matrix. To introduce precisions that give information about the singular vectors, we make the ansatz
\begin{align}
&p(\mathbf{X} | \boldsymbol{\alpha}_L,\boldsymbol{\alpha}_R) \nonumber \\
&= \frac{|\boldsymbol{\alpha}_L|^{q/2} |\boldsymbol{\alpha}_R|^{p/2}}{(2\pi)^{pq/2}} \exp \left( -\frac{1}{2} tr \left( \boldsymbol{\alpha}_L \mathbf{X} \boldsymbol{\alpha}_R \mathbf{X}^\top \right) \right) , \label{gaussianapproximation}
\end{align}
where $\boldsymbol{\alpha}_L \in \mathbb{R}^{p \times p}$ and $\boldsymbol{\alpha}_R \in \mathbb{R}^{q \times q}$ are left and right precision matrices that are positive definite. We use two separate precisions since the left and right singular vectors of a matrix are not equal in general. The ansatz \eqref{gaussianapproximation} is equivalent to setting
\begin{align*}
\mathrm{vec}(\mathbf{X}) \sim \mathcal{N}(\mathbf{0},\boldsymbol{\alpha}_R^{-1} \otimes \boldsymbol{\alpha}_L^{-1}) . 
\end{align*}
To make the precisions promote low-rank, we use Wishart distributions \cite{Bishop} as priors for the precisions, i.e.
\begin{align}
&p(\boldsymbol{\alpha}_L) \propto |\mathbf{S}_L|^{\nu_L/2} |\boldsymbol{\alpha}_L|^{(\nu_L-p-1)/2} e^{-\frac{1}{2} tr(\mathbf{S}_L \boldsymbol{\alpha}_L)} , \label{alphaLprior}\\
&p(\boldsymbol{\alpha}_R) \propto |\mathbf{S}_R|^{\nu_R/2} |\boldsymbol{\alpha}_R|^{(\nu_R-q-1)/2} e^{-\frac{1}{2} tr(\mathbf{S}_R \boldsymbol{\alpha}_R)} , \label{alphaRprior}
\end{align}
where $\mathbf{S}_L,\mathbf{S}_R$ are symmetric and positive definite, $\nu_L > p-1$, $\nu_R > q-1$ and $| \cdot |$ denotes the determinant. The Wishart distributions reduce to gamma distributions if we restrict the precision matrices to be diagonal. Since the Wishart distribution can be seen as a generalization of the gamma-distribution, this choice of priors naturally generalizes the standard RVM. By marginalizing over $\boldsymbol{\alpha}_L$ we get that
\begin{align*}
p(\mathbf{X}|\boldsymbol{\alpha}_R) \propto \frac{1}{\left| \mathbf{X} \boldsymbol{\alpha}_R \mathbf{X}^\top + \mathbf{S}_L \right|^{\nu_L/2}}
\end{align*}
The MAP estimate of $\mathbf{X}$ and the precisions can thus be interpreted as a weighted log-det heuristic.

We model the noise as
\begin{align}
&p(\mathbf{n} | \beta) = \mathcal{N}(0,\beta^{-1} \mathbf{I}_m) , \nonumber \\
&p(\beta) \propto \beta^c e^{-d\beta}, \label{betaprior}
\end{align}
i.e., the noise is modeled as white Gaussian with a gamma distributed precision \cite{Bishop}. We refer to the RVM with the priors \eqref{alphaLprior}, \eqref{alphaRprior} and \eqref{betaprior} as the Relevance Singular Vector Machine (RSVM).

The RSVM iteratively updates the estimate $\mathbf{\hat{X}}$ and the precisions. For fixed precisions, the MAP estimate of $\mathbf{X}$ becomes
\begin{align}
&\mathrm{vec}(\mathbf{\hat{X}}) = \beta \boldsymbol{\Sigma} \mathbf{A}^\top \mathbf{y} ,  \nonumber\\
&\boldsymbol{\Sigma} = ((\boldsymbol{\alpha}_R \otimes \boldsymbol{\alpha}_L) + \beta \mathbf{A^\top A})^{-1} . \label{sigma}
\end{align}


The precisions are updated by maximizing the marginal distribution $p(\mathbf{y},\boldsymbol{\alpha}_L, \boldsymbol{\alpha}_R,  \beta)$. This gives us the update equations
\begin{align}
&\boldsymbol{\alpha}_L^{new} = \nu_L' \left(\boldsymbol{\Sigma}_R + \mathbf{\hat{X}}\boldsymbol{\alpha}_R \mathbf{\hat{X}}^\top + \mathbf{S}_L \right)^{-1} \label{alphaL_update}\\
&\boldsymbol{\alpha}_R^{new} = \nu_R' \left(\boldsymbol{\Sigma}_L + \mathbf{\hat{X}}^\top \boldsymbol{\alpha}_L \mathbf{\hat{X}} + \mathbf{S}_R \right)^{-1} \label{alphaR_update}\\
&\beta^{new} = \frac{m+2c}{||\mathbf{y-A\hat{x}}||_2^2 + tr(\mathbf{A} \boldsymbol{\Sigma}\mathbf{A}^\top) + 2d} \label{beta_update}
\end{align}
where $\nu_L' = \nu_L +q-p-1$, $\nu_R' = \nu_R+p-q-1$ and $\boldsymbol{\Sigma}_R \in \mathbb{R}^{p \times p}$ and $\boldsymbol{\Sigma}_L \in \mathbb{R}^{q\times q}$ are matrices with components
\begin{align*}
&[\boldsymbol{\Sigma}_R]_{kl} = tr(\boldsymbol{\Sigma} (\boldsymbol{\alpha}_R \otimes \mathbf{E}_{kl}^{(L)}))\\
&[\boldsymbol{\Sigma}_L]_{kl} = tr(\boldsymbol{\Sigma} (\mathbf{E}_{kl}^{(R)} \otimes \boldsymbol{\alpha}_L))
\end{align*}
where $\mathbf{E}_{kl}^{(L)} \in \mathbb{R}^{p \times p}$ ($\mathbf{E}_{kl}^{(R)} \in \mathbb{R}^{q \times q}$) is a matrix with $1$ in position $(k,l)$ and zeros otherwise.

We note that the MAP estimate \eqref{sigma} resembles the iterative reweighted least squares approach \cite{Fornasier}, but that the weights are chosen using the Bayesian framework.

One problem with using two precisions is that they can become unbalanced, i.e. one precision can become large and the other small. To balance the precisions, we rescale them as
\begin{align*}
\boldsymbol{\alpha}_L \to \boldsymbol{\alpha}_L \cdot g \cdot h , \,\,\,\, \boldsymbol{\alpha}_R \to \boldsymbol{\alpha}_R \cdot g \cdot h^{-1} , 
\end{align*}
in each iteration, where
\begin{align*}
g = \frac{\sqrt{tr(\boldsymbol{\alpha}_L^{-1}) tr(\boldsymbol{\alpha}_R^{-1})}}{|| \mathbf{\hat{X}}||_F} ,\,\, 
h = \sqrt{||\boldsymbol{\alpha}_R||_F / ||\boldsymbol{\alpha}_L||_F} . 
\end{align*}

The balancing of the precisions removes the dependence on $\nu_L'$ and $\nu_R'$. To make the prior non-informative we take the limits
\begin{align*}
\mathbf{S}_L \to \mathbf{0} , \mathbf{S}_R \to \mathbf{0} . 
\end{align*}
In computations, however, the parameters are given small but non-zero values to avoid numerical instabilities.


\subsection{Symmetric $\mathbf{X}$}

When $\mathbf{X} \in \mathbb{R}^{p \times p}$ is symmetric, the left and right singular vectors are equal (up to sign changes). The left and right precisions can therefore be chosen to be equal. The prior of $\mathbf{X}$ then becomes
\begin{align}
&p(\mathbf{X} | \boldsymbol{\alpha}) = \frac{|\boldsymbol{\alpha}|^{p} }{(2\pi)^{p^2/2}} \exp \left( -\frac{1}{2} tr \left( \boldsymbol{\alpha} \mathbf{X} \boldsymbol{\alpha} \mathbf{X} \right) \right) , \label{Xprior_symmetric}
\end{align}
where $\boldsymbol{\alpha} \in \mathbb{R}^{p \times p}$ has the prior distribution
\begin{align*}
p(\boldsymbol{\alpha}) \propto |\mathbf{S}|^{\nu/2} |\boldsymbol{\alpha}|^{(\nu-p-1)/2} e^{-\frac{1}{2} tr(\mathbf{S} \boldsymbol{\alpha})}
\end{align*}
and $\mathbf{S} \in \mathbb{R}^{p \times p}$ is positive definite. Let
\begin{align*}
\boldsymbol{\Sigma} = \sum_{k=1}^s (\boldsymbol{\Sigma}_{L,k} \otimes \boldsymbol{\Sigma}_{R,k}) ,
\end{align*}
be the Kronecker product decomposition of
\begin{align*}
\boldsymbol{\Sigma} = \left( (\boldsymbol{\alpha} \otimes \boldsymbol{\alpha}) + \beta \mathbf{A^\top A} \right)^{-1}  ,
\end{align*}
where $\boldsymbol{\Sigma}_{L,k}, \boldsymbol{\Sigma}_{R,k} \in \mathbb{R}^{p \times p}$. The update equation for the precision then becomes
\begin{align*}
\boldsymbol{\alpha}^{new} = (\nu-p-1) \left( 2 \mathbf{\hat{X}} \boldsymbol{\alpha}\mathbf{\hat{X}} + \boldsymbol{\Sigma}_R+ \boldsymbol{\Sigma}_L + \mathbf{S} \right)^{-1}  , 
\end{align*}
where
\begin{align*}
\boldsymbol{\Sigma}_L = \sum_{k=1}^s tr(\boldsymbol{\Sigma}_{L,k}\boldsymbol{\alpha}) \boldsymbol{\Sigma}_{R,k} ,\,\,
\boldsymbol{\Sigma}_R = \sum_{k=1}^s tr(\boldsymbol{\Sigma}_{R,k}\boldsymbol{\alpha}) \boldsymbol{\Sigma}_{L,k} . 
\end{align*}

\section{Accelerated RSVM}

We note that the computational complexity of the RSVM is dominated by the computation of $\mathbf{\hat{X}}$ and $\boldsymbol{\Sigma}$ which requires the inversion of a $pq \times pq$ matrix. The computational complexity is thus $\mathcal{O}((pq)^3)$. Using the Woodbury matrix identity \cite{golub96} to rewrite the inverse reduces the complexity to $\mathcal{O}(m^3)$, which can still be large. Here we further reduce the complexity by using the variational Bayesian framework.

Let $\Omega \subset [p] \times [q]$ be a set of indices with corresponding parameters $\mathbf{x}_\Omega$ we want to estimate and let $\Omega^c$ be the complement of $\Omega$.
 We decompose the problem as
\begin{align*}
&\mathbf{A} \mathrm{vec}(\mathbf{X}) = \mathbf{A}_\Omega \mathbf{x}_\Omega +\mathbf{A}_{\Omega^c} \mathbf{x}_{\Omega^c}\\
&tr(\boldsymbol{\alpha}_L \mathbf{X}\boldsymbol{\alpha}_R \mathbf{X}^\top) =\\
& \left( \begin{array}{c}
\mathbf{x}_\Omega \\ \mathbf{x}_{\Omega^c}
\end{array} \right)^\top \left( \begin{array}{cc}
\boldsymbol{\alpha}_{\Omega,\Omega} & \boldsymbol{\alpha}_{\Omega,\Omega^c}\\
\boldsymbol{\alpha}_{\Omega^c, \Omega} & \boldsymbol{\alpha}_{\Omega^c,\Omega^c}
\end{array} \right)\left( \begin{array}{c}
\mathbf{x}_\Omega \\ \mathbf{x}_{\Omega^c}
\end{array} \right) . 
\end{align*}

The variational Bayesian framework gives us that
\begin{align*}
&\log p(\mathbf{x}_\Omega, \mathbf{y} | \boldsymbol{\alpha}_L , \boldsymbol{\alpha}_R) = E_{\mathbf{x}_{\Omega^c}} [ \log p(\mathbf{x}_\Omega, \mathbf{x}_{\Omega^c}, \mathbf{y}| \boldsymbol{\alpha}_L , \boldsymbol{\alpha}_R) ] \\
&= \text{constant} - \frac{\beta}{2} \left|\left| \mathbf{y} - \mathbf{A}_{\Omega^c} (E[\mathbf{x}_{\Omega^c}]) - \mathbf{A}_\Omega (\mathbf{x}_\Omega) \right|\right|_2^2 \\
&- \frac{1}{2} \mathbf{x}_\Omega^\top \boldsymbol{\alpha}_{\Omega,\Omega} \mathbf{x}_\Omega - \mathbf{x}_\Omega^\top \boldsymbol{\alpha}_{\Omega,\Omega^c} E[\mathbf{x}_{\Omega^c}]
\end{align*}
where the constant contains terms which do not depend on $\mathbf{x}_\Omega$.
 By assuming that $\mathbf{x}_{\Omega^c}$ is a random variable with mean $\mathbf{\hat{x}}_{\Omega^c}$ and that it is independent from $\mathbf{x}_\Omega$, we find that the MAP estimate of $\mathbf{x}_\Omega$ becomes
\begin{align}
&\mathbf{\hat{x}}_\Omega = \beta \boldsymbol{\Sigma}_\Omega  \mathbf{A}_\Omega^\top \left( \mathbf{y} - \mathbf{A}_{\Omega^c}  \mathbf{\hat{x}}_{\Omega^c}\right)  - \boldsymbol{\Sigma}_\Omega \boldsymbol{\alpha}_{\Omega,\Omega^c} \mathbf{\hat{x}}_{\Omega^c}  ,  \label{iterative_estimator}\\
&\boldsymbol{\Sigma}_\Omega = \left( \boldsymbol{\alpha}_{\Omega,\Omega} + \beta \mathbf{A}_\Omega^\top \mathbf{A}_\Omega \right)^{-1} . \nonumber
\end{align}
We see that the variational Bayesian approach becomes the block descent method, i.e. the method of minimizing the objective (the negative log likelihood) over different blocks of variables iteratively. Since the objective is convex, the iterative method converges to the minimum of the objective (the full MAP estimate) \cite{Luenberger}.

To accelerate the computation of $\boldsymbol{\Sigma}_R$ and $\boldsymbol{\Sigma}_L$, we make the approximation that the parameter estimates are unbiased and uncorrelated, i.e.
\begin{align*}
& E[(\mathbf{x}_{\Omega_i} - \mathbf{\hat{x}}_{\Omega_i})(\mathbf{x}_{\Omega_j} - \mathbf{\hat{x}}_{\Omega_j})^\top] = \left\{
\begin{array}{ll}
\mathbf{0} & \text{ if } i \neq j\\
\boldsymbol{\Sigma}_{\Omega_i} & \text{ if } i = j
\end{array} \right. , 
\end{align*}
for different blocks $\Omega_i$ and $\Omega_j$ such that $\Omega_i \cap \Omega_j = \emptyset$. This reduces the complexity of calculating $\mathbf{\Sigma}$ from $\mathcal{O}((pq)^3)$ to $\mathcal{O}( \max_i |\Omega_i|^3)$. Since the matrix $\boldsymbol{\Sigma}$ becomes sparse, the computation of $\boldsymbol{\Sigma}_R$ and $\boldsymbol{\Sigma}_L$ can be performed efficiently. The precisions are then updated using \eqref{alphaL_update}, \eqref{alphaR_update} and \eqref{beta_update}.

If we split $\mathbf{X} \in \mathbb{R}^{p \times q}$ into $s \cdot t$ blocks of size $\frac{p}{s} \times \frac{q}{t}$ and the estimator \eqref{iterative_estimator} is iterated over the blocks $K$ times, the number of floating point operations in each iteration reduce from $\mathcal{O}((pq)^3)$ to
\begin{align*}
\mathcal{O}\left(\max \left\{ K \frac{(pq)^3}{(st)^2} , p^3, q^3 \right\}  \right) . 
\end{align*}
However, the reduction in computational complexity comes at the cost of possible loss of accuracy.

\section{Experimental results}
\label{sec:majhead}

To numerically compare the performance of the different algorithms we generated low-rank matrices $\mathbf{X = LR^\top}$ where $\mathbf{L} \in \mathbf{R}^{p \times r}$, $\mathbf{R} \in \mathbb{R}^{q \times r}$ have $\mathcal{N}(0,1)$ components (so $\mathrm{rank}(\mathbf{X}) = r$ with probability $1$).

We measure the performance of the algorithms in terms of the Normalized Mean Square Error
\begin{align*}
\text{NMSE} = \frac{E[||\mathbf{X - \hat{X}}||_F^2]}{E[||\mathbf{X}||_F^2]},
\end{align*}
which was evaluated through numerical simulations. The Signal-to-Noise Ratio (SNR)
\begin{align*}
\mathrm{SNR} = \frac{E[||\mathbf{A} \mathrm{vec} (\mathbf{X})||_2^2]}{E[||\mathbf{n}||_2^2} = \left\{ \begin{array}{lr}
r/ \sigma_n^2 & \text{(completion)}\\
 pqr/ m \sigma_n^2 & \text{(reconstruction)}
\end{array} \right.
\end{align*}
was set to 20 dB. For each scenario we generated $10$ low-rank matrices and $10$ measurements of every matrix for each parameter value.

We set the precisions of the RSVM to identity matrices in the first iteration. For the variational Bayes procedure (VB) we set $s = \min (p,q)$ and for the accelerated RSVM, we partitioned the column vectors of $\mathbf{X}$ into $4$ blocks. For the nuclear norm heuristic \eqref{nuclearnorm} we assumed the noise variance $\sigma_n^2$ to be known and used $\delta = \sigma_n \sqrt{m + \sqrt{8m}}$ as proposed in \cite{Candes10}, where $m$ is the number of measurements. We used the cvx toolbox \cite{cvx} for the nuclear norm minimization.


\subsection{Matrix completion}

\begin{figure}[t]
\includegraphics[width = \columnwidth]{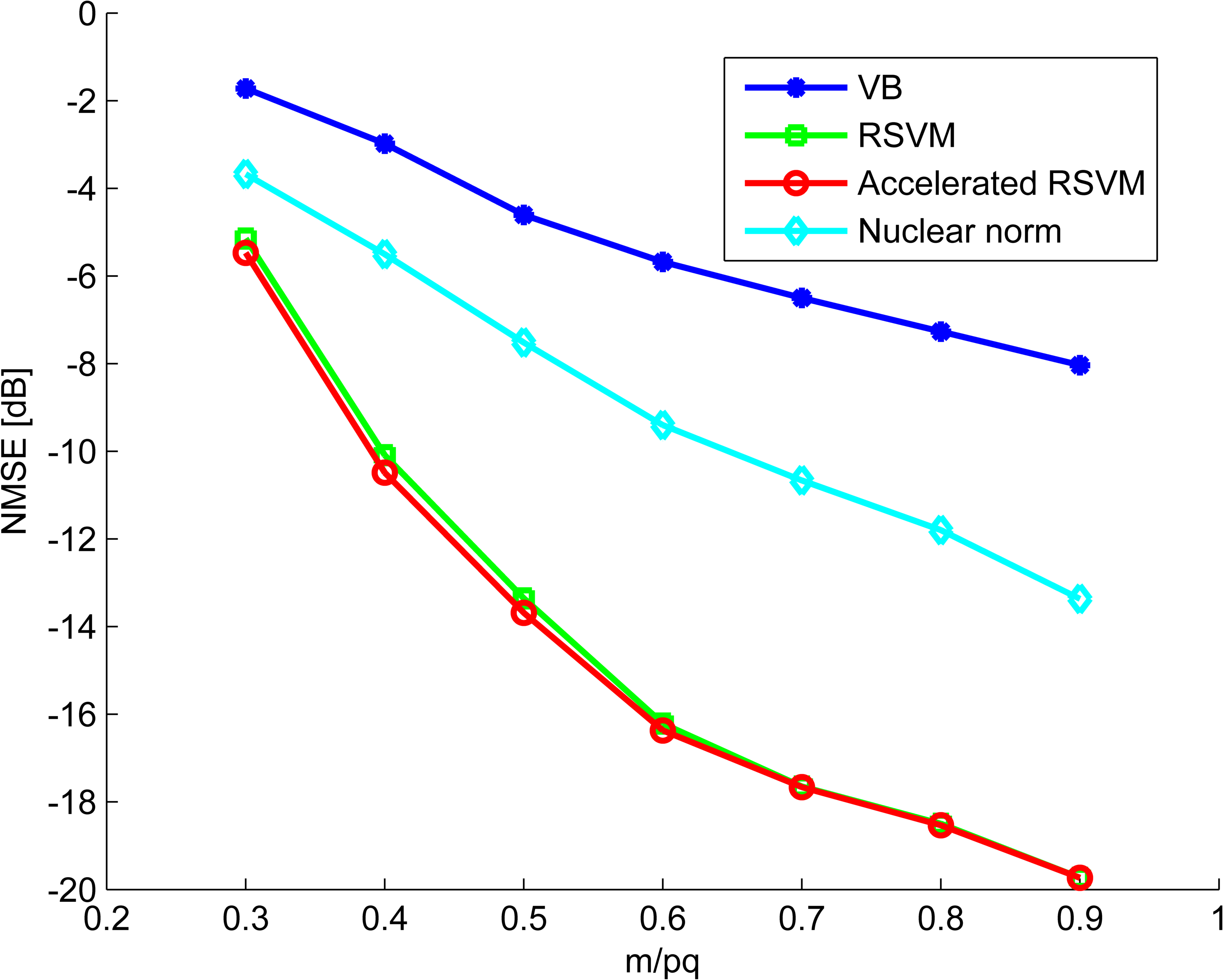}
\caption{Matrix completion. NMSE vs. $m/pq$ for $p=15$, $q = 30$ and $r = 3$.}
\label{completion_rho}
\end{figure}

\begin{figure}[t]
\includegraphics[width = \columnwidth]{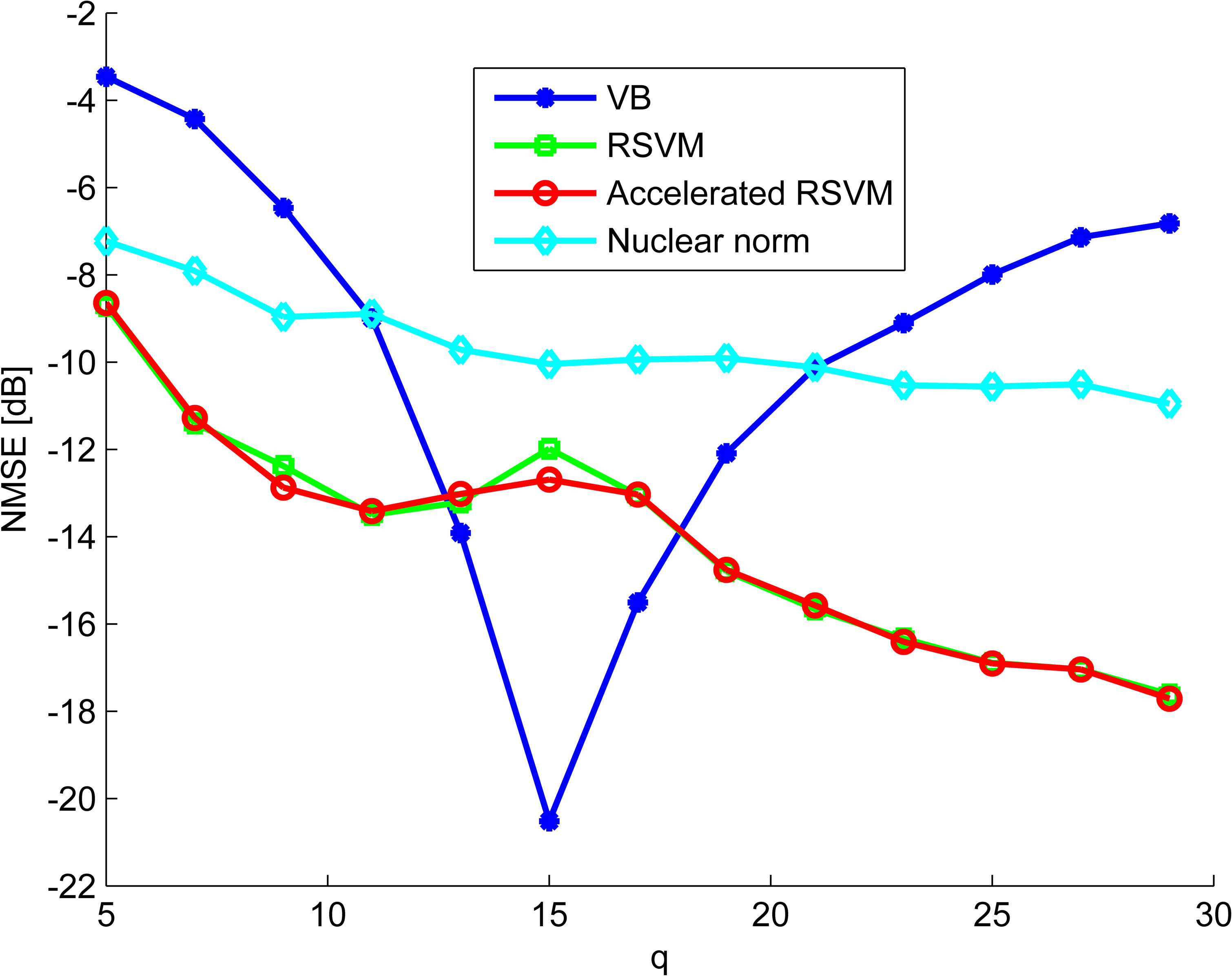}
\caption{Matrix completion. NMSE vs. $q$ for $p=15$, $r=3$, $m/pq = 0.7$.}
\label{completion_q}
\end{figure}


To test the algorithms for the matrix completion problem \cite{Candes09,Candes10,SVT,admira}, $m$ components was chosen uniformly at random from a matrix and the matrix was reconstructed using RSVM, the accelerated RSVM, the variational Bayes procedure (VB) \cite{Babacan} and the nuclear norm \eqref{nuclearnorm}. In the experiments we set $p = 15$, $q = 30$, $r = 3$ and varied the number of measurements $m$.We found that RSVM had $4.7-11.2$ dB lower NMSE than VB and $3.6-6.3$ dB lower NMSE than the nuclear norm for $m/pq \geq 0.4$ (see figure \ref{completion_rho}).


We found that the performance of the variational Bayesian technique in \cite{Babacan} was dependent on ratio of the matrix dimensions. To compare how the performance of the algorithms depend on the ratio we set $p=15$, $r=3$, $m/pq = 0.7$ and varied $q$. We found that VB works better than the other methods when $p=q$, but its performance degrades when $p \neq q$ (see figure \ref{completion_q}). The performance of the nuclear norm, RSVM and the accelerated RSVM is less sensitive to the ratio of the matrix dimensions.


\subsection{Matrix reconstruction}

\begin{figure}[t]
\includegraphics[width = \columnwidth]{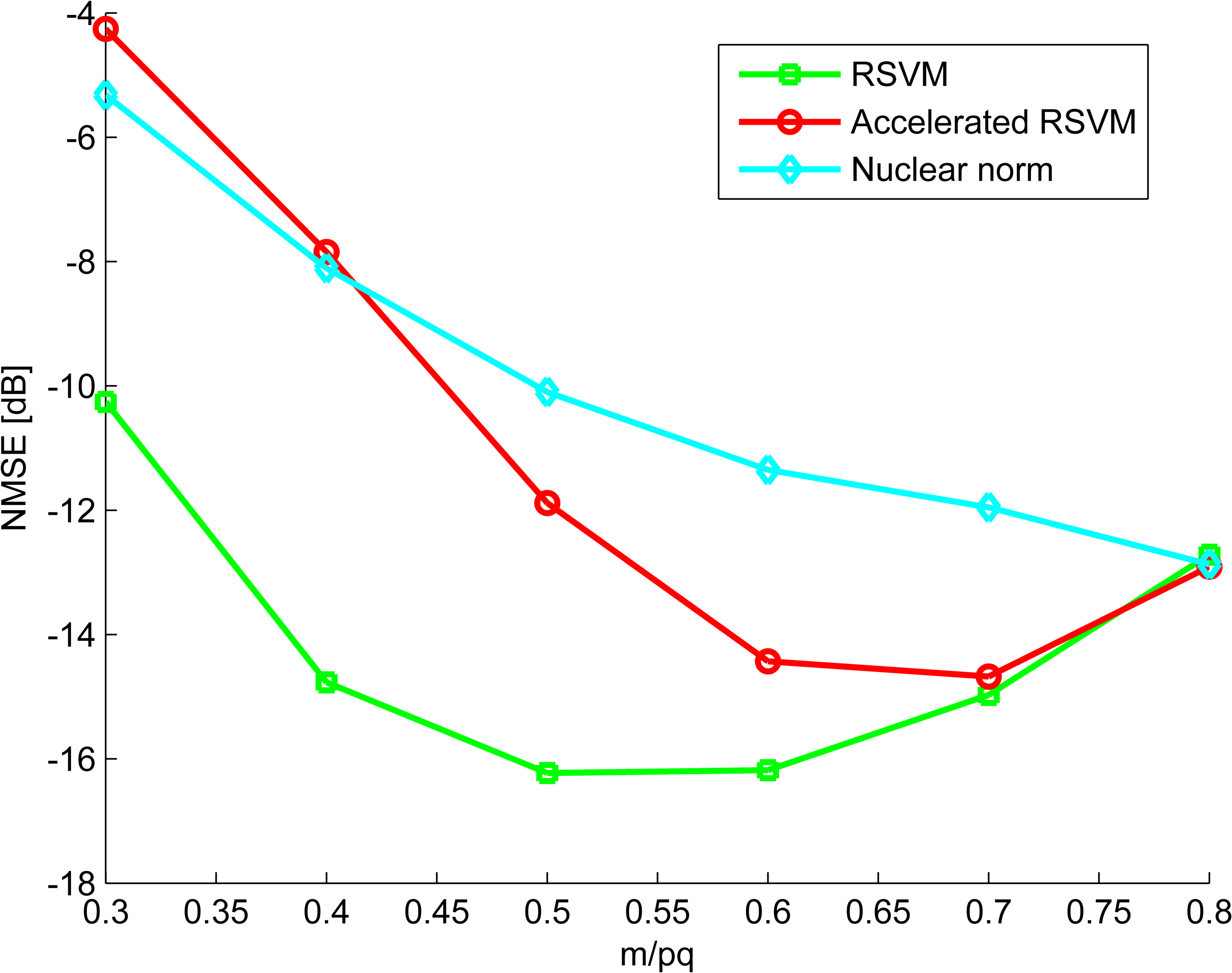}
\caption{Matrix reconstruction. NMSE vs. $\rho$ for $p=q=15$ and $r = 2$.}
\label{rec_rho}
\end{figure}

In matrix reconstruction, the sensing matrix $\mathbf{A}$ in \eqref{measurements} takes linear combinations of the elements in $\mathbf{X}$. 
To examine the performance of the algorithms for this scenario we generated $\mathbf{A} \in \mathbb{R}^{m \times pq}$ by drawing its elements from a $\mathcal{N}(0,1)$ distribution and normalizing the column vectors to unit length. In the simulation we choose $p = q = 15$ and $r = 2$. To estimate $\mathbf{X}$ we used RSVM, the accelerated RSVM and the nuclear norm. We found that RSVM had $5$ dB lower NMSE than the nuclear norm for $0.3 \leq m/pq \leq 0.6$ while the accelerated RSVM had higher NMSE than the nuclear norm for $m/pq \leq 0.4$ and smaller NMSE than the nuclear norm for $m/pq > 0.4$ (see figure \ref{rec_rho}).

\begin{figure}[t]
\includegraphics[width = \columnwidth]{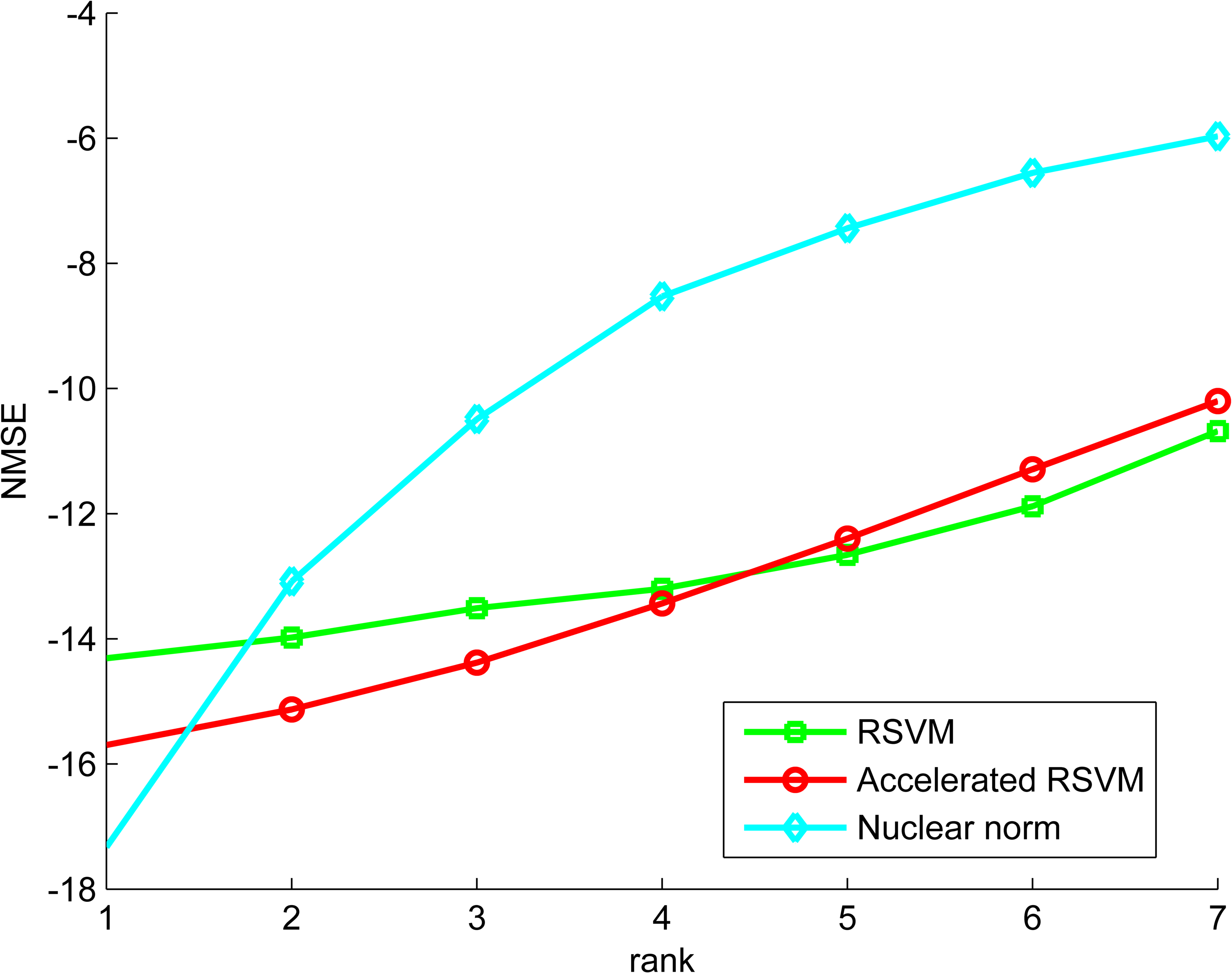}
\caption{Matrix reconstruction. NMSE vs. $\mathrm{rank}(\mathbf{X})$ for $p = q = 15$, $m/pq = 0.7$.}
\label{rec_rank}
\end{figure}

In a second experiment we set $p = q = 15$, $m/pq = 0.7$ and varied the rank of $\mathbf{X}$. We found that the RSVM and the accelerated RSVM had $3-5$ dB lower NMSE than the nuclear norm when $\mathrm{rank}(\mathbf{X}) \geq 3$ (see figure \ref{rec_rank}).

\section{Conclusion}

In the gamut of Bayesian sparse kernel machines, we show that it is possible to use priors on singular vectors to promote low rank. The priors were used to construct a Relevance Vector Machine for low-rank matrix reconstruction. We call the algorithm the Relevance Singular Vector Machine (RSVM). We show that iterative algorithms can be designed to achieve the maximum-a-priori (MAP) solution. Even computational complexity can be further reduced by appropriate approximations leading to accelerated algorithms.

Through simulations, it was shown that the RSVM outperformed the variational Bayesian method for matrix completion when the matrix is skewed. The RSVM also outperformed the relaxed nuclear norm heuristic in the matrix reconstruction scenario.


\end{document}